 \UseRawInputEncoding
\documentclass[11pt,hyp]{article}
\RequirePackage{doi}
\usepackage{hyperref}

\usepackage{varwidth}

\usepackage[leqno]{amsmath}
\usepackage{amssymb,amsthm,graphicx,enumitem,upref}

\usepackage[titletoc, toc, title]{appendix}
\usepackage{breakcites}
\numberwithin{equation}{section}
\theoremstyle{plain}
\newtheorem{thm}{Theorem}[section]

\theoremstyle{definition} 
\newtheorem{exam}[thm]{Example}

\usepackage{lipsum}

\usepackage{enumitem}
\newlist{myQuoteEnumerate}{enumerate}{1}% Set max nesting depth
\setlist[myQuoteEnumerate,1]{label=(2.\arabic*)}% Use numbers for level 1

\newenvironment{MyQuote}{%
    \begin{myQuoteEnumerate}[resume=*,series=MyQuoteSeries]%
    \item \begin{quote}%
}{%
    \end{quote}%
    \end{myQuoteEnumerate}%
}%

\begin{document}

\title{Fundamental Errors in Kane and Mertz's Alleged Debunking\\ 
of Greater Male Variability in Mathematics Performance
}

\author{Rosalind Arden and Theodore P. Hill}

\date{\vspace{-5ex}}  

\maketitle

\begin{abstract}
Kane and Mertz's 2012 \textit{AMS Notices} article ``Debunking Myths about Gender and Mathematics Performance" claims to have debunked the greater male variability hypothesis with respect to mathematics abilities. The logical and statistical arguments supporting their claim, however, which are being widely cited in the scientific literature, contain fundamental errors. The methodology is critically flawed, the main logical premise is false, and the article omits reference to numerous published scientific research articles that contradict its findings. Most critically, Kane and Mertz's final conclusion that their data are inconsistent with the greater male variability hypothesis is wrong. The goal of the present note is to correct the scientific record with respect to those claims. Most importantly, by publicizing these errors, the \textit{Notices} will reduce the chance of similar future errors being repeated.
\end{abstract}

\section{Introduction}

The American Mathematical Society's sole article of record addressing the controversial ``hypothesis" of greater male variability (GMV) is the widely cited paper ``Debunking Myths About Gender and Mathematics Performance" by husband and wife team Jonathan Kane and Janet Mertz \cite{vob68}.

That feature-length article which appeared in the \textit{Notices of the AMS}, ``the world's most widely read magazine aimed at professional mathematicians" seems to promise, through its sensational title, a dramatic trouncing of false collective beliefs regarding gender gaps in mathematics performance. Yet surprisingly the body of the article makes no mention of myths, but rather presents an investigation of at least five \textit{hypotheses} (italics theirs) based on data from  children's performance on standard tests. Kane and Mertz conclude that most of these hypotheses are ``inconsistent" with their findings, thus apparently ``debunking" these pernicious myths. The initial spotlight is focused on greater male variability, suggesting it is the first and foremost ``myth" to be exposed.

Unfortunately, the Kane-Mertz article is seriously flawed in virtually every aspect of the GMV discussion.  
The goal of the present note is \textit{not} to argue the validity or invalidity of GMV in general or with respect to cognitive or mathematical abilities in humans.  Rather it is  to correct the scientific record with respect to the basic errors in their article. Perhaps this might also inspire others to take a closer look at the rest of their hypothesis/myths, especially given their concern about how resources are spent in these various endeavors.
Most importantly, by publicizing these errors, the \textit{Notices} will reduce the chance of similar future errors being repeated.

\section{Misleading History}
The lead paragraph of \cite{vob68} introduces the GMV hypothesis as follows:
%\begin{quotation}\noindent
\begin{MyQuote}
``the \textit{greater male variability hypothesis}, originally proposed by Ellis in 1894 and reiterated in 2005 by Lawrence Summers when he was president of Harvard University, states that variability in intellectual abilities is intrinsically greater among males" \cite[p.~10, italics in original]{vob68}. 
\end{MyQuote}
%\end{quotation}

This statement (2.1) is misleading in both its history and its formulation. First, the GMV hypothesis dates back to Charles Darwin, a fact  explicitly acknowledged by variability expert Stephanie Shields in \cite{vob67}, the only historical article cited in \cite{vob68}. And second, the classical GMV hypothesis pertains to traits in many species throughout the animal kingdom, not only to mathematics performance or even general cognitive abilities in humans. As also noted in \cite{vob67}, the working definition of the hypothesis has changed back and forth over the years, but the authors of this note know of no other source where GMV is interpreted exclusively in terms of variability of mathematics ability or performance. 

By only mentioning ``intellectual abilities", Kane and Mertz explicitly misrepresent  
Havelock Ellis, who gathered the scientific literature on human sex differences in variability ``in order to study the issue directly and in-depth" \cite[p.~775]{vob67}. In Ellis's classic text he devoted an entire chapter to ``Variational Aspects of Men", in which he clearly states ``Both the physical and the mental characters of men show wider limits of variation than do the physical and mental characters of women" \cite[p.~358]{vob11}.

Similarly, by referring only to intellectual variability, Kane and Mertz also misrepresent Harvard President Larry Summers, who in fact said 
%\begin{quotation}\noindent
\begin{MyQuote}
``It does appear that on many, many different human attributes -- height, weight, propensity for criminality, overall IQ, mathematical ability, scientific ability -- there is relatively clear evidence that whatever the difference in means -- which can be debated -- there is a difference in the standard deviation, and variability of a male and a female population" \cite[p.~3]{vob16}.
\end{MyQuote}
%\end{quotation}

The standard form of the GMV hypothesis, as recorded by Darwin, is simply that throughout the animal kingdom, males are generally more variable than females for many traits. It does not say that there is greater male variability in every trait of every species. In the special case where the species is human, the GMV hypothesis simply says that for many traits, both physical and cognitive, male variability is greater than female variability. 

\section{Flawed Methodology}

Kane and Mertz claim they ``tested the greater male variance hypothesis" with respect to mathematics performance in humans, specifically referring to Fields medalists and ``empowerment as reflected by percentage of women in technical, management, and government positions" \cite[p.~12]{vob68}.  That is, their experiments were designed to draw conclusions about \textit{adult} humans, not infants or children.

In the description of their method, however, Kane and Mertz clearly state that ``most measures of mathematics performance here are based on the TIMSS, a quadrennial study that includes a mathematics assessment" of fourth and eighth graders from numerous countries \cite[p.~11]{vob68}.  Kane and Mertz's conclusions are thus mostly based on tests of pre- and early-adolescent \textit{children}, not adults, and it has been established that at those ages boys and girls follow different developmental trajectories in many traits such as school performance  \cite{vob201} and height \cite{vob191}.  

For example, Arden and Plomin \cite{vob17} addressed questions about over- and under-representation of each gender at the low and high extremes of measures of cognitive abilities by studying sex differences in variance of test scores across childhood. Among other conclusions, they found that ``From age 2 to age 4, girls in our study were highly significantly over-represented in the top tail" \cite[p.~44]{vob17}. Thus, employing the same logic and methodology as \cite{vob68} to extrapolate data from tests done on children from age 2 to 4 to conclusions about adults, this finding of Arden and Plomin would imply that among adults, \textit{women} are highly significantly over-represented in the top tail of intelligence. Similarly, extrapolating predicted heights of adults from data about the fourth and eighth graders that \cite{vob68} studied would conclude that adult women are taller on average than men. To draw reasonable inferences about gender differences in variability (or any other traits) among human adults from tests on children requires serious formal justification, and this is missing in \cite{vob68}.

 Similarly, the working definition of the GMV hypothesis explicitly stated in  the Kane and Mertz article is that ``variability \textit{in intellectual abilities} is intrinsically greater among males"  \cite[p.~12, emphasis added]{vob68}, yet as noted above their conclusions are based solely on standard tests of  \textit{mathematical ability}. To infer that mathematical ability alone is a reasonable measure of overall intrinsic  intellectual ability is highly debatable and requires justification. This, too, is entirely missing in  \cite{vob68}.

\section{Crucial Logical Error}

To quantify the notion of gender differences in variability of a collection of data, the standard parameter used is the so-called \textit{variance ratio} VR, which, by convention, is defined as the variance of the male data divided by the variance of the female data \cite[p.~11]{vob68}.
There is greater male variability in a particular trait of a given sexually-dimorphic species if the VR for that trait is greater than one. In Kane and Mertz's main logical argument, they assert 

\begin{quotation}\noindent
``we tested the greater male variance hypothesis. If true, the variance ratios (VRs) for all countries should be greater than unity and similar in value" \cite[p.~13]{vob68}.
\end{quotation}

This logical argument is repeated almost verbatim in a ``Doing the Math" section of the American Academy of Arts and Sciences flagship publication \textit{Science} article highlighting their work: ``If the greater male variability hypothesis \ldots is true, then that variability would persist, consistently, across all 86 countries" \cite[p.~2]{vob166}.   In short, the main logical premise of \cite{vob68} is this: 

\begin{center}
\begin{varwidth}{\textwidth}
\begin{enumerate}[label = (\alph*)]
\item[(P)\hspace{1.5em}] If there is GMV for humans worldwide,  i.e., if \textit{VR} $>1$ for a particular trait,\\
                                        then the corresponding VRs for each country should all be greater than 1 and similar in value.                                                                            
\end{enumerate}
\end{varwidth}
\end{center}

Although (P) may appear intuitive and plausible, it is simply false. The premise (P)  violates a standard statistical fact, namely, the formula for the variance of a finite weighted mixture of distributions, as will be seen in the next two examples. \\

The first example is hypothetical and illustrates how a union of countries (in this case only two) could exhibit greater male variability as a whole, even though not all of the individual countries do. That is, $\mathit{VR} >1$ for the union, but not for each country.  
\begin{exam}
\label{ex2}
A population consists of two countries $C_1$ and $C_2$, with equal numbers of people in each, divided equally in each among men and women.  Scores on a certain test administered to everyone in the overall population result in means $m_1$ and $m_2$ and standard deviations $\sigma_1$ and $\sigma_2$ for the men in countries $C_1$ and $C_2$, respectively, and means $f_1$ and $f_2$ and standard deviations $\hat{\sigma}_1$ and $\hat{\sigma}_2$ for women.  Applying the standard formula (e.g., equation (1.21) in \cite{novob32}) for the moments of finite weighted mixtures of distributions, the variance of the men's scores in the overall population is given by $({2(\sigma_1^2 + \sigma_2^2) + (m_1 - m_2)^2)}/4$ and that of women is $({2(\hat{\sigma}_1^2 + \hat{\sigma}_2^2) + (f_1 - f_2)^2)}/4$.
Letting $\mathit{VR}$ denote the variance ratio in the overall population and $\mathit{VR}_i$ the variance ratios in $C_i$, $i=1,2$, it follows immediately that 

\begin{equation}
\label{eq1}
\mathit{VR_1} = \frac{\sigma_1^2}{\hat\sigma_1^2}, 
\mathit{VR_2} =  \frac{\sigma_2^2}{\hat\sigma_2^2},
\mathit{VR} = \frac{2(\sigma_1^2 + \sigma_2^2) + (m_1 - m_2)^2}{2(\hat\sigma_1^2 + \hat\sigma_2^2) + (f_1 - f_2)^2}.
\end{equation}

For example, if $m_1 = m_2 = 102, \sigma_1^2 = 5, \sigma_2^2 = 1, f_1 = f_2 = 103, \hat\sigma_1^2 = 1, \hat\sigma_2^2 = 2$, then 
equation (\ref{eq1}) implies that $\mathit{VR_1} = 5,  \mathit{VR_2} = 0.5$, and $ \mathit{VR} = 2$, so there is greater male variability in the overall population, but greater female variability in $C_2$. This contradicts (P).
\end{exam}

The second example is also based on humans, but in contrast to the previous example, uses real data and concerns the physical trait of height. (Human height is one of science's most studied and documented measurements and has been recorded and analyzed in great detail, over time and geographic location, in part because height is easy to measure and is an indicator of important factors such as nutrition and genetics.) 

\begin{exam}
Roser, Appel and Ritchie \cite{vob191} list the mean height of men worldwide at 178.4 cm with a standard deviation of 7.59 cm and women's mean height at 164.7 cm with a standard deviation of 7.07 cm.  The variance ratio for adult human height worldwide is therefore $\mathit{VR}>1.07$, implying greater male variability.  The variance ratios by country and birth year, on the other hand, are ``all over the place" and range from less than 0.5 to greater than 2.5 \cite{novob24, novob23}.  Thus there is greater variability worldwide in heights of men than heights of women, even though the VR's for height vary significantly among countries.  This also contradicts (P).
\end{exam} 

This false premise (P) is being propagated in the scientific literature. For example, a 2021 article in the \textit{Journal of Comparative Economics} specifically cites Kane and Mertz \cite{vob68}, and employs their faulty logic:
\begin{quotation}
\noindent
``the `male greater variability' hypothesis does not accommodate the staggering cross-country differences found here. Kane and Mertz (2012), examining mathematics performance, also find that the male-to-female performance variance ratio significantly differs across countries, which is inconsistent with the greater male variability hypothesis"  \cite[p.~438]{vob204}.
\end{quotation}
  
No justification for (P) is presented in \cite{vob68}, and as just seen in the two previous counter-examples,  (P) \textit{is false prima facie}. 

\section{Erroneous Main Conclusion}
A logical argument can be wrong, yet its conclusion correct, but this is not the case in \cite{vob68}. The authors themselves explicitly concede that they found 8\% greater male variability in their analysis of mathematics performance as indicated by certain test scores worldwide - ``These findings agree well with the VR of 1.08 reported from a large meta-analysis involving data from 242 studies involving over 1 million Americans" \cite[p.~13]{vob68}. 
This finding was confirmed by Reilly \textit{et al} who reported that their own study of global gender differences in variance in mathematics and science achievement was consistent with the greater male variability hypothesis, and that ``Similar patterns were observed by Kane and Mertz (2012)" \cite[p.~42]{vob129}. 

In their conclusions, however, Kane and Mertz maintain exactly the opposite: ``These findings are inconsistent with the greater male variability hypothesis" \cite[p.~14]{vob68}.   This conclusion was repeated in an even stronger form in  \textit{Science}, where Mertz asserted \\

\begin{varwidth}{\textwidth}
\begin{enumerate}[label = (\alph*)]
\item[(C)\hspace{2em}] ``We have pretty clear data debunking the greater male variability hypothesis"  \cite[p.~1]{vob166}.
\end{enumerate}
\end{varwidth}
\\
 
Kane and Mertz also conclude that the non-uniformity in VR's they found is largely an artifact of  ``a complex variety of sociocultural factors rather than intrinsic differences" \cite[p.~11]{vob68}, i.e., cultural factors as opposed to ``innate, biologically determined differences between the sexes" \cite[p.~10]{vob68}. 

This same secondary conclusion is repeated in the \textit{AAAS Science} article which reports that 
``cross-cultural analysis seems to rule out several causal candidates, including coeducational schools, low standards of living, and  innate variability among boys" \cite{vob166}.  (Note, however, that it is only innate variability that is emphatically ``debunked").

Kane and Mertz again repeat this claim in \textit{Scientific American}:
\begin{quotation}
\noindent
``The finding that males' variance exceeds females' in some countries but is less than females' in others and that both range `all over the place suggests it can't be biologically innate, \textit{unless you want to say that human genetics is different in different countries},' Mertz argues.  `The vast majority of the differences between male and female performance must reflect social and cultural  factors.' \cite[p.~4, emphasis added]{vob159}".
\end{quotation}

But evidence shows that  human genetics \textit{do} differ within and among countries, even within the continent of Europe \cite{novob29}.  In Kane and Mertz's study, Taiwan and Tunisia are seen to have the extreme variance ratios in 2007 {TIMSS} scores for eighth graders, namely, 1.31 for Taiwan and 0.91 for Tunisia \cite[Table~2]{vob68}.  How do Kane and Mertz conclude that the significantly different VR's they found for these two countries are primarily artifacts of sociocultural factors rather than, say, a more balanced combination of sociocultural and innate biological factors? This is an empirical question that has been studied in some depth (e.g., \cite{vob120}); it is not a matter of ``what you want to say".
Yet based on \cite{vob68}, \textit{Scientific American} concludes ``Now that  the greater male variability hypothesis has fallen short, nature is not looking as important as scientists once thought" \cite[p.~5]{vob159}.
  
\section{Failure to Cite Published Counter-evidence}

The Kane and Mertz article \cite{vob68} fails to report that their own GMV historical reference clearly states that in 19th century studies ``The biological evidence \textit{overwhelmingly} favored males as the more variable sex" \cite[p.~773, emphasis added]{vob67}. Similarly, their article \cite{vob68}  includes 53 references, but omits scores of research studies that already had reported greater male variability in many different contexts, both cognitive and otherwise. 

Even for the very special case of mathematics performance, Kane and Mertz  fail to report that their same historical reference states that  ``Research from the 1960s and 1970s \ldots [indicates that] scores on tests of mathematical ability show a consistent trend in the direction of greater male variability" \cite[p.~793]{vob67}. As for more recent decades, Kane and Mertz also failed to cite more than a dozen \textit{previously-published} studies supporting GMV in mathematical or quantitative cognitive abilities (emphasis added in the following):
\begin{quotation}
\noindent
``Finally, it should be noted that the boys' \textbf{SAT-M} scores had a larger variance than the girls'" \cite[p.~1031]{vob104}. 
\end{quotation}

\begin{quotation}
\noindent
``The important exception to the rule of vanishing gender differences is that the well-documented gender gap at the upper levels of performance on high school \textbf{mathematics} has remained constant over the past 27 years" \cite[p.~95]{vob193}. 
\end{quotation}

\begin{quotation}
\noindent
``Not only do males reliably score higher on the \textbf{SAT-M} (mean difference approximately .5 standard deviation), they also display a greater variability on such measures" \cite[p.~328]{vob194}. 
\end{quotation}

\begin{quotation}
\noindent
``Males were consistently more variable than females in \textbf{quantitative reasoning, spatial visualization}, spelling, and general knowledge" \cite[p.~61]{vob75}. 
\end{quotation}

\begin{quotation}
\noindent
``The current finding that males were more variable than females in \textbf{math and spatial abilities} in some countries is consistent with the findings of greater male variability in these abilities in the United States" \cite[p.~90]{vob23}. 
\end{quotation}

\begin{quotation}
\noindent
``Examination of the ratios of male score variance to female score variance (VR values) in Table 2 [including \textbf{Mathematics}] reveals that the variance of male scores is larger than that of female scores" \cite[pp.~43--44]{vob27}. 
\end{quotation}

\begin{quotation}
\noindent
``\textbf{As in mathematics}, the variance of the total scores among boys [in science] was generally larger than that among girls across all participating countries" \cite[p.~371]{vob71}. 
\end{quotation}

\begin{quotation}
\noindent
``A nationally representative UK sample of over 320,000 school pupils aged 11--12 years was assessed on \ldots separate nationally standardized tests for verbal, \textbf{quantitative}, and non-verbal reasoning \ldots for all three tests there were substantial sex differences in the standard deviation of scores, with greater variance among boys" \cite[pp.~463,475]{vob37}. 
\end{quotation}

\begin{quotation}
\noindent
``Males are more variable on most measures of \textbf{quantitative and visuospatial ability}, which necessarily results in more males at both high- and low-ability extremes" \cite[p.~1]{vob25}. 
\end{quotation}

\begin{quotation}
\noindent
``Despite the modest differences at the center of the distribution, the greater variability of male scores resulted in large asymmetries at the tails, with males out-numbering females by a ratio of 7 to 1 in the top 1\% on tests of \textbf{mathematics and spatial reasoning} \ldots greater male variance is observed even prior to the onset of preschool" \cite[pp.~220--221]{vob19}. 
\end{quotation}

\begin{quotation}
\noindent
``This indicates that the differences that we observed -- for example, in the overrepresentation of males at the extremes of the distributions for \textbf{quantitative reasoning} \ldots With one exception \ldots all variance ratios were greater than 1.0" \cite[p.~395]{vob31}. 
\end{quotation}

\begin{quotation}
\noindent
``The hypothesis of greater male variability was supported in most domains [including \textbf{Arithmetic}]" \cite[p.~475]{vob132}. 
\end{quotation}

\begin{quotation}
\noindent
``[T]he methods are applied to estimate the residual variance in test scores for male and female students on the \textbf{mathematics portion} of the 2007 Arizona Instrument to Measure Standards \ldots We find that male students exhibit greater residual as well as raw variability for this data set" \cite[pp.~2937, 2947]{vob133}. 
\end{quotation}
   
\begin{figure}[!ht] 
  \center\includegraphics[width=0.8\textwidth]{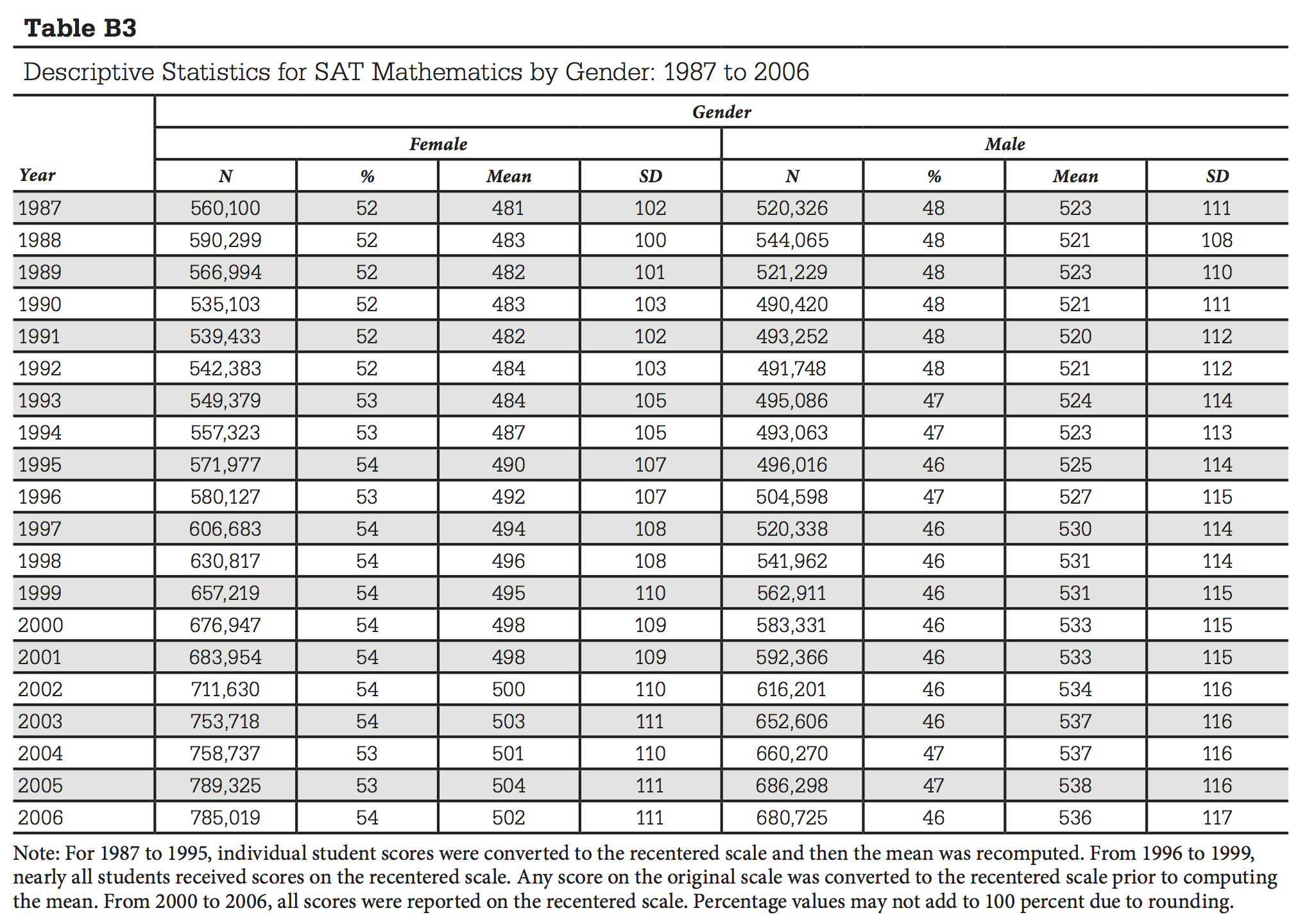}
  \caption{Table B3 from College Board Report 2006-5 \cite[p.~25]{novob30}. For this twenty year period, the VRs are all between 1.09 and 1.22 with average about 1.14.}
  \label{kobrin1}
\end{figure}

Kane and Mertz also omit an important report that had been published by the College Board (see Figure \ref{kobrin1}).  This official historical summary of 20 years of results of SAT-M, the mathematics portion of the Scholastic Aptitude Test, shows \textit{greater male variability in every single year}, with average VR approximately 1.14.  

Note that the male standard deviations are not much larger than those of females, but the mean values of males are also noticeably higher.  This combination can have significant effects on the right tails since, as Feingold pointed out, ``what might appear to be trivial group differences in both variability and central tendency can cumulate to yield very appreciable differences between the groups in numbers of extreme scorers" \cite[p.~11]{vob59}; see also Example 3.6 in  \cite{vob192}. 

\section{Summary}

One of the main objectives of \cite{vob68} is to claim that their data prove that greater male variability (with respect to mathematics ability in particular) is ``bunk" - a myth or a fairy tale which if true, may help explain the preponderance of male Fields medalists.  
Yet greater male variability also leads to more males being found in the left or lower tail of a distribution, thus perhaps contributing to the fact that seventy per cent of the children eligible for special education are boys \cite[p.~97]{vob190}, and that boys and men are overrepresented at the lowest levels of IQ \cite[p.~42]{vob22}, \cite[p.~529]{vob29}.  Greater male variability has implications for \textit{both} tails of the distribution, and this important fact is essentially ignored in \cite{vob68}.

Recall that it is not  the goal of the present note to argue the validity or invalidity of greater male variability in general or with respect to cognitive or mathematical abilities in humans, but simply to correct the scientific record concerning the faulty arguments and conclusions in \cite{vob68}. However, the interested reader should note that since the publication of their paper \cite{vob68} in 2012, overwhelming additional new evidence of GMV in various traits and species has been published; in particular, here are quotes from several very recent mega-studies. The first even uses data from the same PISA tests used by Kane and Mertz:
\begin{quotation}
\noindent
``Twelve databases from IEA [International Association for the Evaluation of Educational Achievement] and PISA [Program for International Student Assessment] were used to analyze gender differences within an international perspective from 1995 to 2015 \ldots The `greater male variability hypothesis' is confirmed" \cite[p.~1]{vob10}.
\end{quotation}

%The second may perhaps help explain the possibility that ``variability in intellectual abilities is intrinsically greater among males" \cite[pp.~10]{vob68} is true, and may have biological underpinnings:

The second suggests a possible biological contribution to greater male variability:

\begin{quotation}
\noindent
 ``the largest-ever mega-analysis of sex differences in variability of brain structure, based on international data spanning nine decades of life \ldots The present study included a large lifespan sample and robustly confirmed previous findings of greater male variance in brain structure in humans. We found greater male variance in all brain measures, including subcortical volumes and regional cortical surface area and thickness, at both the upper and the lower end of the distributions" \cite[pp.~6,23]{vob167}.
\end{quotation}

As emphasized above,  GMV for humans refers to many traits, and strong evidence supporting this phenomenon has also appeared recently. For example, this study corroborates Professor Summers's observation (2.2) above almost exactly: 

\begin{quotation}
\noindent
``The principal finding is that human intrasex variability is significantly higher in males, and consequently constitutes a fundamental sex difference \ldots The data presented here show that human greater male intrasex variability is not limited to intelligence test scores, and suggest that \textit{generally greater intrasex variability among males is a fundamental aspect of the differences between sexes}" \cite[pp.~220--221, emphasis added]{vob5}.
\end{quotation}
(Evidence of GMV in cognitive abilities has also recently been  reported for the first time in non-human  \cite{vob173} and even non-mammalian species \cite{vob174}.) \\

The arguments above pointed out fatal methodological and logical errors in the widely cited Kane and Mertz article \cite{vob68}, and that their main conclusion (C) is false. Their data do \textit{not} debunk the greater male variability hypothesis, even in the special case of mathematics performance. On the contrary, the preponderance of new evidence suggests exactly the opposite and the reader is encouraged to look at Appendix A of \cite{vob4}, which surveys GMV studies published since 2000. The \textit{Notices of the AMS} is looked up to by scientists in many fields, and should be held to especially high standards, including corrections of published and widely propagated errors. The goal of this note is to help correct the scientific record and illuminate  errors involving the greater male variability hypothesis, thereby reducing the likelihood that future similar international studies repeat the same basic errors.

%%%%%%%%%%%%%%%%%%%%%%%%%%%%

\end{document}